\documentclass[utf8]{FrontiersinHarvard} 

\usepackage{url,hyperref,lineno,microtype,subcaption}
\usepackage[onehalfspacing]{setspace}
\usepackage{amsmath,amssymb,amsfonts}
\usepackage{graphicx}
\usepackage{algorithm,algorithmic}
\usepackage[algo2e]{algorithm2e}
\usepackage{hyperref}
\hypersetup{hidelinks=true}
\usepackage{textcomp}
\usepackage{color,colortbl}
\usepackage{soul,xcolor}
\usepackage{multirow}
\usepackage{subcaption}

\definecolor{green1}{rgb}{0.4660 0.6740 0.1880}
\colorlet{trgreen}{green1!30}

\definecolor{blue1}{rgb}{0 0.4470 0.7410}
\colorlet{trblue}{blue1!30}

\definecolor{yellow1}{rgb}{0.9290 0.6940 0.1250}
\colorlet{tryellow}{yellow1!30}

\definecolor{red1}{rgb}{0.8500 0.3250 0.0980}
\colorlet{trred}{red1!30}

\definecolor{purple1}{rgb}{0.4157 0.0510 0.6784}
\colorlet{trpurple}{purple1!30}

\colorlet{trteal}{teal!30}

\DeclareMathOperator*{\argmin}{arg\,min}

\def\keyFont{\fontsize{8}{11}\helveticabold }
\def\firstAuthorLast{Sirlanci {et~al.}} 
\def\Authors{Melike Sirlanci\,$^{1,2,*}$, George Hripcsak\,$^{3}$, Cecilia C. Low Wang\,$^{4}$, J. N. Stroh\,$^{1}$, Yanran Wang\,$^{5}$, Tellen D. Bennett\,$^{1}$, Andrew M. Stuart\,$^{6}$ and David J. Albers\,$^{1,3,5,7}$}
\def\Address{$^{1}$ Department of Biomedical Informatics, School of Medicine, University of Colorado Anschutz Medical Campus, Aurora, CO, USA\\
$^{2}$Department of Applied Mathematics, University of Colorado Boulder, Boulder, CO, USA\\
$^{3}$Department of Biomedical Informatics, Columbia University, New York, NY, USA\\
$^{4}$Division of Endocrinology, Metabolism and Diabetes, Department of Medicine, School of Medicine, University of Colorado Anschutz Medical Campus, Aurora, CO, USA\\
$^{5}$Biostatiscs \& Informatics Department, Colorado School of Public Health, Aurora, Colorado, USA\\
$^{6}$Department of Computing and Mathematical Sciences, California Institute of Technology, Pasadena, CA, USA\\
$^{7}$ Department of Bioengineering, University of Colorado Denver, Denver, CO, USA}
\def\corrAuthor{Melike Sirlanci}
\def\corrEmail{melike.sirlanci@cuanschutz.edu}

\begin{document}
\onecolumn
\firstpage{1}

\title[Personalized Model-Based Glycemic Management]{A Stochastic Model-Based Control Methodology for Glycemic Management in the Intensive Care Unit} 

\author[\firstAuthorLast ]{\Authors} 
\address{} 
\correspondance{} 

\extraAuth{}

\maketitle

\begin{abstract}

\section{}
Intensive care unit (ICU) patients exhibit erratic blood glucose (BG) fluctuations, including hypoglycemic and hyperglycemic episodes, and require exogenous insulin delivery to keep their BG in healthy ranges. Glycemic control via glycemic management (GM) is associated with reduced mortality and morbidity in the ICU, but GM increases the cognitive load on clinicians. The availability of robust, accurate, and actionable clinical decision support (CDS) tools reduces this burden and assists in the decision-making process to improve health outcomes. Clinicians currently follow GM protocol flow charts for patient intravenous insulin delivery rate computations. We present a mechanistic model-based control algorithm that predicts the optimal intravenous insulin rate to keep BG within a target range; the goal is to develop this approach for eventual use within CDS systems. In this control framework, we employed a stochastic model representing BG dynamics in the ICU setting and used the linear quadratic Gaussian control methodology to develop a controller. We designed two experiments, one using virtual (simulated) patients and one using a real-world retrospective dataset. Using these, we evaluate the safety and efficacy of this model-based glycemic control methodology. The presented controller avoids hypoglycemia and hyperglycemia in virtual patients, maintaining BG levels in the target range more consistently than two existing GM protocols. Moreover, this methodology could theoretically prevent a large proportion of hypoglycemic and hyperglycemic events recorded in a real-world retrospective dataset.

\tiny
 \keyFont{ \section{Keywords:} personalized stochastic model, modeling blood glucose dynamics, glycemic management in the intensive care unit, clinical decision support, reducing cognitive burden of healthcare professionals} 
\end{abstract}

\section{Introduction}\label{intro}

Glycemic management (GM) in the intensive care unit (ICU) is a challenging and crucial task. The rapidly changing health condition of patients and frequent interventions result in highly nonstationary blood glucose (BG) behavior, which means that when considered as a time series, their mean, variance, and covariance change over time. Maintaining BG levels in a target range is associated with reduced morbidity and mortality. However, performing GM in ICU is a demanding job for caregivers. Thus, providing efficient clinical decision support (CDS) can potentially improve GM outcomes and reduce the workload of caregivers. ICU patients often receive fluids and medications, such as glucocorticoids, pressors including norepinephrine and epinephrine, fluoroquinolones, and beta-blockers, which affect their insulin sensitivity. All these factors make the GM of ICU patients very challenging. An efficient and actionable model-based GM strategy must address these ICU-specific challenges and limitations.

Our goal is to develop a \textit{personalized model-based glucose control method for eventual use as a CDS tool in the ICU}. Since a multitude of factors which generate the resulting nonstationary BG behavior cannot easily be incorporated into mechanistic models, we aim to model the BG dynamics at a sufficient resolution to provide useful GM strategy in settings where these factors are simply represented as stochastic uncertainty. We developed a linear quadratic Gaussian (LQG) controller based on a linear stochastic BG model, the Minimal Stochastic Glucose (MSG) model \citep{sirlanci2023simple}. This model was developed specifically for the purpose of glycemic management in ICU. Thus, we use the MSG model as the basis of this control methodology because \textit{(i)} it can be used for all patients regardless of their diabetic status, \textit{(ii)} it has an analytical solution that allows use in online settings, \textit{(iii)} it mitigates parameter identifiability issues that occur mostly when an unmeasured system component (e.g., interstitial insulin) is included as a state of the dynamical system model, \textit{(iv)} its stochasticity allows for robust control, and \textit{(v)} its ability to track the mean BG behavior and the amplitude of BG oscillations provides clinical interpretation.

We designed different experiments with simulated and real retrospective data to evaluate this controller. The pipeline can easily be adjusted for different GM strategies, i.e., the LQG controller can be tuned to address different target glycemic regions and account for factors such as nutritional intake. The structure of the simulated data experiment can be used to test the safety and efficacy of different model-based GM strategies, providing a tool for the investigation of new GM strategies. The developed methodology addresses ICU-specific challenges, such as rapidly changing BG behavior and sparse data. In this sense, the methodology is transferable other ICU settings such as hemodynamic management, with appropriate modification of relevant mechanistic models.

Our \textit{contribution} in this paper is summarized below.

\begin{itemize}
	\item We build the control methodology using the LQG control technique based on the previously developed MSG model for personalized GC of ICU patients who require IV insulin delivery and are fed through an enteral tube.
	\item We develop the retrospective data evaluation framework to test the potential effectiveness of the model-based GC methodology using previously collected real-world ICU data.
\end{itemize}

\subsection{Background}

Healthcare professionals follow flow charts, called GM protocols, to maintain patient BG levels within a specific target range. These protocols vary across ICU types (e.g., medical or neurological) and between institutions. The flow chart uses recent BG values and the intravenous (IV) insulin rate at the intervention time to calculate the new IV insulin delivery rate and subsequent BG measurement time. These protocols were developed based on clinical trials designed to regulate ICU patient BG behavior in response to glycemic control (GC) regimens and may expose controversial GM differences in the ICU amongst institutions.

Some clinical trials showed that the intensive insulin therapy (IIT) reduced mortality and/or mortality among ICU patients when compared to conventional insulin therapy (CIT) \citep{van2001intensive,van2006intensive}. Even though these results changed how GM performed in the ICU, these results could not be replicated by follow-up clinical trials \citep{annane2010corticosteroid,arabi2008intensive,brunkhorst2008intensive,coester2010intensive,de2008strict,green2010intensive,macrae2014randomized,preiser2009prospective}. A comprehensive clinical trial, the NICE-SUGAR trial \citep{nice2009intensive}, showed that the IIT increased mortality among adult ICU patients.

Although there is no definitive explanation for these disparate results, potential hypotheses have been offered \citep{clain2015glucose}, the most widely accepted of which is that increased hypoglycemia clouds the benefits of tight glycemic control (TGC). Besides causing dangerous complications such as coma and death, hypoglycemia may also result in irreversible complications such as neuronal damage and cardiac arrhythmia. \textit{Therefore, a crucial feature of a glycemic controller is its ability to avoid hypoglycemia.}

These GM protocols were developed with a one-size-fits-all approach. However, recent research emphasizes the potential of personalized healthcare \citep{mathur2017personalized,nardini2021evolution}. Several studies showed that GM results in significantly different outcomes between different patient groups, such as those with or without diabetes, \citep{green2012hyperlactatemia,kaukonen2014stress,krinsley2013diabetic,lanspa2013moderate,leibowitz2010effects}. These differences in response to the same therapy among different patient groups could be one of the reasons for the contradictory results of TGC in clinical trials. One study \citep{suhaimi2010makes} comparing TGC performance in two different ICUs concluded that TGC protocols should be developed using intra- and inter-patient variability as well as carbohydrate administration. So far, there has been no clinical trial to investigate these differences between patients. Using computational modeling tools to develop glycemic controllers, including baseline evaluation by virtual trials, may reveal BG behavior in different patient populations and aid the development of patient group-specific GC strategies.

The GC problem has been widely investigated in the artificial pancreas (AP) project, which focuses on automatic and semi-automatic insulin delivery via model-based control algorithms in the type 1 diabetes mellitus (T1DM) context \citep{sorensen1985physiologic,hovorka2004nonlinear,dalla2007meal,weinzimer2008fully,hovorka2011closed,ghorbani2014reducing,man2014uva,goharimanesh2015fractional,pinsker2016randomized,delavari2019adaptive,moon2021current}. Translating similar ideas to the ICU requires addressing specific challenges: dynamic BG behavior, insulin secreted by the body independently of external insulin delivery, viz., endogenous insulin, and data sparsity. Perhaps the most significant challenge to GC in the ICU is highly variable BG behavior. This variability is caused by the body's stress response to severe illness, interventions, constant feeding, and exogenous insulin delivery when needed. Even though the factors affecting BG levels, such as corticosteroids, are known, the resulting dynamics are generally not. Therefore, these factors are not incorporated into physiology-based mechanistic models of the glucose-insulin system in the ICU, which increases model error. Moreover, even if the mechanistic relationship between these factors and BG levels is known, inclusion in the model increases complexity and causes computational intractability.

Another data-related limitation is sparse BG measurement data (e.g., $\sim$10-15 measurements per day). The highly fluctuating nature of BG levels means that $\sim$10-15 measurements are not enough to accurately track patient glycemic dynamics, and continuous glucose monitors are not routine in the ICU. Therefore, an effective model-based glycemic controller for ICU patients must: (a) \textit{be safe, i.e., avoid hypoglycemia and hyperglycemia}; (b) \textit{account for inter- and intra-patient variability}; and (c) \textit{address challenges specific to the ICU setting}. In every step of the glycemic controller development process presented here, we aim to address all these challenges.

\subsection{Literature Review}

Researchers have investigated GM in the ICU setting from various perspectives. This literature review is organized according to different aspects of GC in the ICU, and concludes with some observations about the use of ideas from control in personalized medicine.

Many research groups designed clinical trials to test the safety and eﬀicacy of IIT compared to CIT. Depending on the specific ICU setting, the target BG ranges can vary. Some of these trials reported benefits of IIT over CIT \citep{van2001intensive,van2006intensive} and some resulted in the opposite conclusions or failed to confirm previous results \citep{arabi2008intensive,brunkhorst2008intensive,de2008strict,nice2009intensive,preiser2009prospective,coester2010intensive,green2010intensive,macrae2014randomized}. Also, some studies suggested that moderate glycemic control (MGC) could be more beneficial for all or certain patient groups than TGC \citep{clain2015glucose,lanspa2013moderate,leibowitz2010effects}.

None of these GM protocols account for nutrition or for inter- and intra-patient variability. However, researchers claimed these protocols should account for these types of variability among patients, \citep{suhaimi2010makes,chase2011physiological}. In \citep{chase2007model,lonergan2006simple}, the authors developed a GM protocol that regulates both the exogenous insulin and nutrition delivery. They
117 tested the efficacy and feasibility of this protocol in a pilot study \citep{lonergan2006pilot}.

Some researchers developed physilogical mechanistic models \citep{van2006minimal,hovorka2008simulation,haverbeke2008nonlinear,pielmeier2010simulation,lin2011physiological,zhou2023model}, and control-theoretic algorithms based on these models to optimize and personalize the IV insulin rate for more effective GM, \citep{chase2018improving}. In \citep{stewart2016safety, knopp2021goldilocks}, the authors compared the eﬀicacy of GM protocols using retrospective data. Virtual or \textit{in silico} clinical trials provide a means to test the safety and eﬀicacy of these algorithms before pilot studies and clinical trials. Several studies found \textit{in silico} testing useful, advantageous and safe to test the efficacy and feasibility of GM approaches and for validation purposes \citep{wilinska2008silico,wilinska2011evaluating,chase2010validation,fisk2012star,stewart2018nutrition,uyttendaele2020risk}.

Some of these control algorithms have been further tested in pilot studies \citep{pielmeier2010glucosafe,gonzalez2022effectiveness} and in clinical trials \citep{hovorka2007blood,pachler2008tight,pielmeier2012decision,van2013logic,dubois2017software,uyttendaele2019star,uyttendaele2021star,uyttendaele2023clinical}. In all these studies, the algorithm-based GC approaches provided safe GM and resulted in improved GC evaluated by several different measures including time spent in target range, mean BG values, and number of adverse events. These results show evidence that personalized treatment could increase the success rate at the individual level compared to standard one-fits-all approaches.

Finally, a recent research program studies the statistical foundations of personalized approaches to medicine \citep{liang2023randomization}; building on this it includes study of a personalized control-theoretic approach to drug delivery, in the specific context of how to diminish delivery of addictive drugs \citep{gradu2023online}. Similar ideas are relevant in our work.

\subsection{Outline}

Section \ref{data} describes the retrospective ICU dataset; Section \ref{comp_method} provides the details of our computational methodology; Section \ref{exp_design} describes the experiments designed for simulated and real-world data settings, and Section \ref{evaluation} presents the evaluation metrics used to test safety and eﬀicacy of the model-based controller. Then, Section \ref{num_res} presents the results. Finally, Section \ref{discussion} discusses the strengths, advantages, and limitations of this personalized model-based GC methodology and Section \ref{conc} summarizes our findings.

\section{Materials and Methods}\label{data_exp}

\subsection{Data}\label{data}

The data were extracted from the University of Colorado Health Data Compass Clinical Data Warehouse between 2010 and 2019 and represent six units, including medical, burn, surgical/trauma, neurological, cardiothoracic, and cardiac ICUs. The inclusion criteria for this study were patients on a tube-feed who were in the ICU for at least three days and whose interval between the first and last recorded IV insulin delivery was at least one day. If a particular patient has two intervals that satisfy these criteria, broken up by the removal and insertion of the enteral tube, then we treat these data as two separate experiment time intervals (and similarly for multiple such intervals). We excluded patients who were pregnant. Notably, we did not exclude anyone based on their diabetes status. These data resulted in a collection of 106 experiment time intervals drawn from 106 patients since none of the patients had multiple ICU stays that met the inclusion criteria. Note that the insulin delivered via IV route is typically short-acting. All of these 106 patients were delivered short-acting IV insulin. For evaluation purposes, we also focused on a subset of these data, restricted to patients with at least one hypoglycemic or hyperglycemic episode. These data included 126 adverse events, including 19 hypoglycemic and 107 hyperglycemic events. These refined data represented 23 time intervals belonging to 23 different patients. Detailed information about the patient cohort can be found in Table \ref{data_table}. We also provide this information for all 106 patients in Supplementary Material.  This study was approved by the Colorado Multiple Institutional Review Board with protocol number 18-2519 on January 10, 2023.

\begin{table}[!ht]
	\centering
	\caption{Demographic and ICU stay-related health record information of patients included in the real-world dataset.}
	\begin{tabular}{| p{0.45\linewidth} | r| r | r | r | r |}
		\hline
		& \textbf{mean} & \textbf{stdev} & \textbf{min} & \textbf{median} & \textbf{max}\\ \hline
		\textbf{Age (years)}	& 56.4 		&	13.5				&	25.0			&	60.0					& 78.0\\ \hline
		\textbf{Length of ICU stays that contain an interval meeting the inclusion criteria (days)}	& 12.5	&	11.6	&	4.8	&	8.8	&	51.7\\ \hline
		\textbf{Number of BG measurements}	&	142.6	&	94.4	&	48.0	&	114.0	&	462\\ \hline
		\textbf{Length of tube feed administration (hours)}	& 163.0	&	90.1	&	73.7	&	126.3.0	&	447.0\\ \hline
		\textbf{Length of IV insulin administration (hours)}	& 140.3	&	90.5	&	33.2	&	115.1	&	444.0\\ \hline
		\textbf{Race}		& \multicolumn{5}{l|}{American Indian and Alaska Native: 1}\\
		& \multicolumn{5}{l|}{Asian: 1}\\
		& \multicolumn{5}{l|}{Black or African American: 3}\\
		& \multicolumn{5}{l|}{White or Caucasian: 18}\\ \hline
		\textbf{Ethnicity}		& \multicolumn{5}{l|}{Hispanic: 2}\\
		& \multicolumn{5}{l|}{Non-Hispanic: 21}\\ \hline
		\textbf{Sex}	& \multicolumn{5}{l|}{Female: 10}\\
		& \multicolumn{5}{l|}{Male: 13}\\ \hline
	\end{tabular}
	\label{data_table}
\end{table}

\subsection{Computational Methodology}\label{comp_method}

We parse the control problem into two sub-problems: \textit{(i) the estimation of the unknown model parameters in a personalized way based on the available data (model identification)}, and \textit{(ii) the estimation of the personalized optimal controller based on the identified model}. In the language of control theory these correspond to system identification and optimal control. We will establish a methodology combining several algorithms to solve each sub-problem efficiently. We use the MSG model, introduced in \citep{sirlanci2023simple}, to represent BG dynamics. We use an optimization approach to identify the unknown model parameters and use LQG control method to estimate the optimal IV insulin rate to keep BG level in the desired range.

The MSG model used to represent the BG dynamics is a variant of the well-known Ornstein-Uhlenbeck stochastic process and is analytically solvable. We use the discrete-time form of the model for compatibility with the available data; in this form, the output becomes a multivariate normally distributed random variable. The model consists of two components: a deterministic and a stochastic component. The deterministic component describes the body's own effort to reach the basal glucose level and incorporates the effects of nutrition and exogenous insulin on the BG level. This component aims to model the dynamics that can be resolvable with the available data, e.g., it does not aim to model the detail of BG oscillations. However, since BG levels are known to oscillate, the stochastic component of the model quantifies the magnitude of BG oscillations. This is merely a technical choice that allows us to resolve the dynamics at a level enough for BG forecasting and control with sparsely available data.

\subsubsection{Identification Problem}

We estimate unknown model parameters using patient-specific data to obtain a personalized BG model. Let $\{y_i\}_{i=0}^K$ and $\{x_i\}_{i=0}^K$ represent the BG measurement data and the model output, respectively, collected at time points $\{t_i\}_{i=0}^K$. Denoting data and the model output as vectors, $y=[y_i]_{i=0}^K$ and $x=[x_i]_{i=0}^K$, we have $x\sim \mathcal{N}(\mu,\Sigma)$ where $\mu\in\mathbb{R}^{K+1}$, $\Sigma\in\mathbb{R}^{{K+1}\times{K+1}}$, and $\mathcal{N}$ refers to normal distribution. Note that $\mu$ and $\Sigma$ are functions of the unknown model parameters, $\theta\in\Theta$ where $\Theta$ is the set of admissible parameters. Then, we formulate the parameter estimation problem as follows:
\begin{equation*}
	\theta^* = \argmin_{\theta\in\Theta}\frac{1}{\sqrt{\det{\Sigma(\theta)}}}\exp\left(-\frac{1}{2}(y-\mu(\theta))^T \Sigma(\theta)^{-1} (y-\mu(\theta)) \right).\label{opt_problem}
\end{equation*}
This is a specific example of the more general system identification problems, \citep{aastrom2012introduction,aastrom2021feedback,bechhoefer2021control}. We use MATLAB's \texttt{fmincon} function to solve this constrained optimization problem. More details can be found in \citep{sirlanci2023simple}.

\subsubsection{Control Problem}

We estimate the optimal controller, i.e., the input we can adjust, to maintain desirable system behavior. Several control-theoretical approaches can be used in this setting. To start with the most straightforward method that meets our requirements and to exploit the fact that the MSG model is linear and stochastic, we use an LQG controller \citep{anderson2007optimal,murray2009optimization,aastrom2012introduction}. Consider the linear input/output system
\begin{equation}
	\begin{aligned}
		\dot{x} &= Ax + Bu + \xi, \ \ \xi\sim N(0,C_{\xi}),\\
		y &= Cx + \eta, \ \  \eta\sim N(0, C_{\eta}),\label{control_sys}
	\end{aligned}
\end{equation}
where $A\in \mathbb{R}^{n\times n}$, $B\in \mathbb{R}^{n\times p}$, $C\in \mathbb{R}^{q\times n}$. $\xi$ represents the disturbances to the system ($C_{\xi}\in\mathbb{R}^{n\times n}$) and $\eta$ represents the measurement noise ($C_{\eta}\in\mathbb{R}^{q\times q}$), and they are assumed to be uncorrelated. In this system, $x$ is the state that we aim to ``control" by adjusting the variable $u$ appropriately based on the observed $y$. In our case, $n,p,q=1$ because the BG is the only model state, the exogenous insulin is a single-valued function, and the only collected measurements are BG values.

The LQG controller is the combination of the linear quadratic regulator (LQR) and linear quadratic estimator (LQE or Kalman filter) \citep{murray2009optimization}. The optimal controller for system \eqref{control_sys} has the form
\begin{align*}
	\dot{\hat{x}} &= A\hat{x} + Bu + K_f(y - C\hat{x}),\\
	u &= -K_c(\hat{x} - x_r) + u_r,
\end{align*}
where $x_r$ is the reference value (the target BG value), and $u_r$ is the shift in the controller (IV insulin input) to ensure convergence to the reference value. The \textit{separation principle} for this formulation states that $K_f$ is the optimal observer gain (or Kalman gain) ignoring the controller, and $K_c$ is the optimal controller gain ignoring the noise \citep{anderson2007optimal,aastrom2012introduction}. This principle allows us to perform the state and controller estimation consecutively and incorporate the result of one into the other in the computation process. We provide more details about the algorithmic structure of this process in Section \ref{exp_design}.

\begin{figure}[htb]
	\centering
	\includegraphics[width=\linewidth]{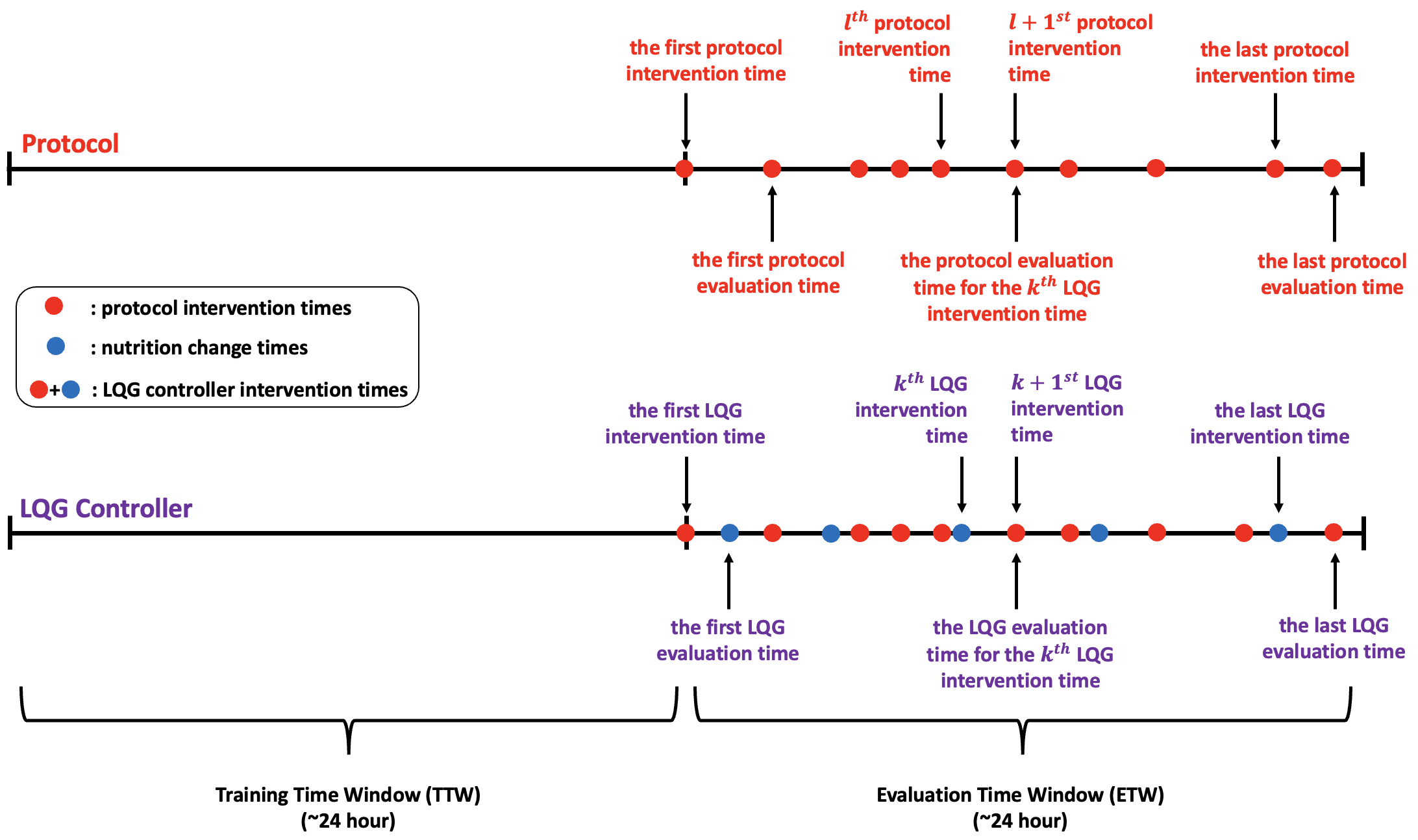}
	\caption{These timelines show the protocol, LQG controller intervention times, and corresponding evaluation times for each intervention time. The LQG controller intervention times are the union of the protocol intervention times and the nutrition change times over the ETW.}
	\label{sim_data_schema}
\end{figure}

\subsection{Design of the Control Framework}\label{exp_design}

We evaluate the eﬀicacy of our BG control strategy by using both simulated and real-world data. This strategy allows us to investigate different aspects of this control methodology. The simulated data experiment is essential for two reasons. First, we can simulate any number of patients for robust evaluation. Second, we can observe and fairly compare the implication(s) of any IV insulin rate by simulating the resulting BG level. The real data experiment allows us to evaluate how the developed control strategy performs with the limitations of real-world data, including data sparsity, measurement noise and rapidly changing BG behavior. Also, using retrospective data required significant changes to the algorithmic structure and evaluation techniques. This is because unlike the simulated data case, we cannot ``deliver" the predicted amount of insulin to observe its effect. So, we have to devise a method to evaluate the controller's performance based on delivered IV insulin to the patients.

We use two different GM protocols in both the simulated data and the real-world data cases. One of them is used in the neurological ICU at the New York - Presbyterian Hospital (NYPH) and can be found in \citep{hripcsak2022evaluating}. The other protocol is based on the Yale insulin infusion protocol and used in several ICUs at the UCHealth University of Colorado Hospital (UCH). These two protocols have different target ranges, 140 - 180 mg/dL and 120 - 150 mg/dL. Our aim in using these protocols is to evaluate the performance of our model-based GM methodology. Since we do not aim to compare the effectiveness of these two protocols, we will not name them in presenting our results. We will call them \textit{Protocol A} and \textit{Protocol B} without identifying which one is which.

\subsubsection{Simulated Data Experiment}\label{sim_data}

We used the ICU Minimal Model (ICUMM) introduced in \citep{vanherpe07,vanherpe06} to generate virtual patient data that have the same elements as real-world data (tube-fed nutrition rate, IV insulin rate, and BG measurements). Using a different model from the MSG model, which is the basis of the control algorithm, is crucial to evaluate controller eﬀicacy. Using the same model to simulate patient data and for the control algorithm would provide a perfect representation of actual BG dynamics with very accurate BG control, but would not mimic real-world conditions or allow for a realistic evaluation of the control algorithm.

For each patient, data from the first 24 hours, the \textit{training time window} (TTW), is used for model estimation (generate the personalized model). We define the \textit{intervention times} as the times at which the IV insulin rate is determined using either the protocol or the machinery developed here. The \textit{evaluation times} are the BG measurement times occurring after an intervention. We use these BG measurements at the evaluation times to assess ``decision" eﬀiciency (i.e., the IV insulin rate delivery at the intervention times). We define the time window encompassing these evaluation time points as the \textit{evaluation time window} (ETW). The intervention times for the LQG controller are the collection of protocol intervention times and nutrition change times over the ETW to account for the nutrition changes in our control approach. A visual description of this experimental design can be seen in Fig. \ref{sim_data_schema}.

\begin{algorithm}[H]
	Generate M sets of nutrition data covering both the TTW (training time window) and ETW (evaluation time window);\\
	Generate M sets of IV insulin data covering only TTW;\\
	\For{each patient}{
		Create patient profile: Generate the ICUMM model parameters by selecting uniformly over their physiological ranges. These parameters characterize each virtual patient\;
		
		Generate BG measurement times over TTW by using the estimated distribution of the time interval lengths between BG measurements from a real-world BG measurement data collected from ICU patients\;
		
		By using the generated nutrition, IV insulin data and model parameters, simulate BG measurements at the predefined BG measurement time points (in the previous step) by using ICUMM\;
	}
	\caption{\textbf{Data generation steps for the simulated data experiment.}}
	\label{data_gen_alg}
\end{algorithm}

\paragraph{Data Generation} Let N denote the total number of virtual patients and M denote the number of different simulated nutrition and IV insulin data. We select M to be smaller than N to mitigate the effect of simulated nutrition and IV insulin data on the results. Then, we randomly create M groups of patients, each group containing N/M patients who are delivered the same nutrition over both the TTW and ETW and IV insulin over only TTW. First, we simulated the nutrition and IV insulin rates using estimated distribution of these inputs from a real-world dataset. To generate these data, we estimated the distributions of the lengths of time intervals over which the nutrition and insulin rates were constant and the distributions of delivered nutrition and insulin rates. Based on these estimated distributions, we generated M different nutrition and IV insulin delivery data. Then we create the patient profiles as follows. We generate each ICUMM model parameter using a uniform distribution over the given feasible range for each parameter. Then, we pair these parameters to create N sets of random model parameters, representing N different virtual patient profiles. Next, we estimated the distribution of the length of time between consecutive BG measurement times and generated BG measurement times for each patient based on this distribution. Finally, we randomly paired the virtual patients and the nutrition-insulin data and simulated the resulting BG values at the generated BG measurement times over the TTW for each patient (with ICUMM). This process produced the simulated data ready to apply control algorithms.

\paragraph{Performing Control} We start with applying the GM protocol to each patient. At each intervention time, we identify and deliver the protocol-suggested IV insulin rate to the virtual patient by simulating their BG value at the next protocol-suggested intervention time (using ICUMM). Note that this is an evaluation time for the most recent intervention time. We perform this process over the whole ETW, recording all the intervention times at this step to use for the LQG control.

For the LQG control, we first define the set of intervention times. As mentioned above, it is important to incorporate the effect of nutrition into the control algorithm. The nutrition rate at the time of intervention is naturally included in the computation of the optimal IV insulin rate via the LQG controller. If a change in the nutrition rate occurred between two protocol-suggested intervention times, we added that nutrition-change time as an intervention time for the LQG control approach.

\begin{algorithm}[H]
	\For{each patient}{
		\textbf{Protocol/Flow Chart:} Let current intervention time = the end point of TTW (training time window)\;
		
		\While{current intervention time $\leq$ the end point of the ETW (evaluation time window)}{
			Use flow chart to compute the IV insulin rate ($I^{protocol}$) and the next BG measurement time, which is also the next intervention time\;
			
			Simulate BG over [current intervention time, next intervention time] by ICUMM to obtain the BG value at the next evaluation time\;
			
			Set current intervention time = next intervention time\;
		}
		\textbf{LQG Control:} Here the intervention times are the union of the intervention times obtained in the protocol step and the nutrition-change times\;
		
		\For{each intervention time}{
			Use the prior 24-hour data to estimate unknown MSG model parameters\;
			
			Use LQG controller with the current BG value and the estimated model parameters to estimate optimal IV insulin rate ($I^{lqg}$)\;
			
			Simulate BG over [current intervention time, next intervention time] by ICUMM to obtain BG value at the next evaluation time\;
		}
	}
	\caption{\textbf{GM steps via protocol/flow chart and model-based controller for the simulated data experiment.}}
	\label{sim_data_control_alg}
\end{algorithm}

The LQG control algorithm requires specifying a target value to estimate the optimal IV insulin rate. To avoid hypoglycemic episodes, we set this value to be the upper bound of the respective protocol (150 mg/dL and 180 mg/dL). In addition, another precaution we took to avoid hypoglycemia was to set a threshold for the LQG controller suggested IV insulin rate. If the controller suggested a rate higher than 25 U/hr, the optimal IV insulin rate was set to 25 U/hr, where U represents the unit of insulin.

At each intervention time, we use the most recent 24-hour data to obtain the personalized MSG model. Then, using this model with the LQG controller, we estimate the optimal IV insulin rate and deliver this insulin rate to the virtual patient by simulating their BG value at the next intervention time (using ICUMM). We perform this process over the entire ETW. The process is described in Algorithms \ref{data_gen_alg} and \ref{sim_data_control_alg}.

\subsubsection{Real Data Experiment}\label{real_data_setting}

In this case, we use retrospective real-world data, described in Section \ref{data}, to investigate the applicability and evaluate the model-based controller's performance in a real-world setting. We designed this experiment by addressing the real-world retrospective data limitations, e.g., intervention times are defined by data, which are the recorded IV insulin change times. Also, we evaluate the performance of the GM techniques by comparing the suggested IV insulin rates and the actual IV insulin rate at the time of interventions, which caused a hypoglycemic or hyperglycemic event.
$ $

\begin{algorithm}[H]
	\For{each patient}{
		Define intervention times: For each hypoglycemic and hyperglycemic measurement, find the latest IV insulin rate change time\;
		
		\textbf{Protocol/Flow Chart:} 
		\For{each intervention time}{
			Use the flow chart to compute the IV insulin rate ($I^{protocol}$) suggested by the protocol\;
		}
		\textbf{LQG Control:}
		\For{each intervention time}{
			Use prior 24-hour data to estimate unknown MSG model parameters\;
			
			Use LQG controller with the current BG value and the estimated model parameters to estimate optimal IV insulin rate ($I^{lqg}$)\;
		}
	}
	\caption{\textbf{GM steps with the protocol/flow chart and the model-based controller using retrospective data.}}
	\label{real_data_cont_alg}
\end{algorithm}

\begin{figure}[ht]
	\centering
	\includegraphics[width=\linewidth]{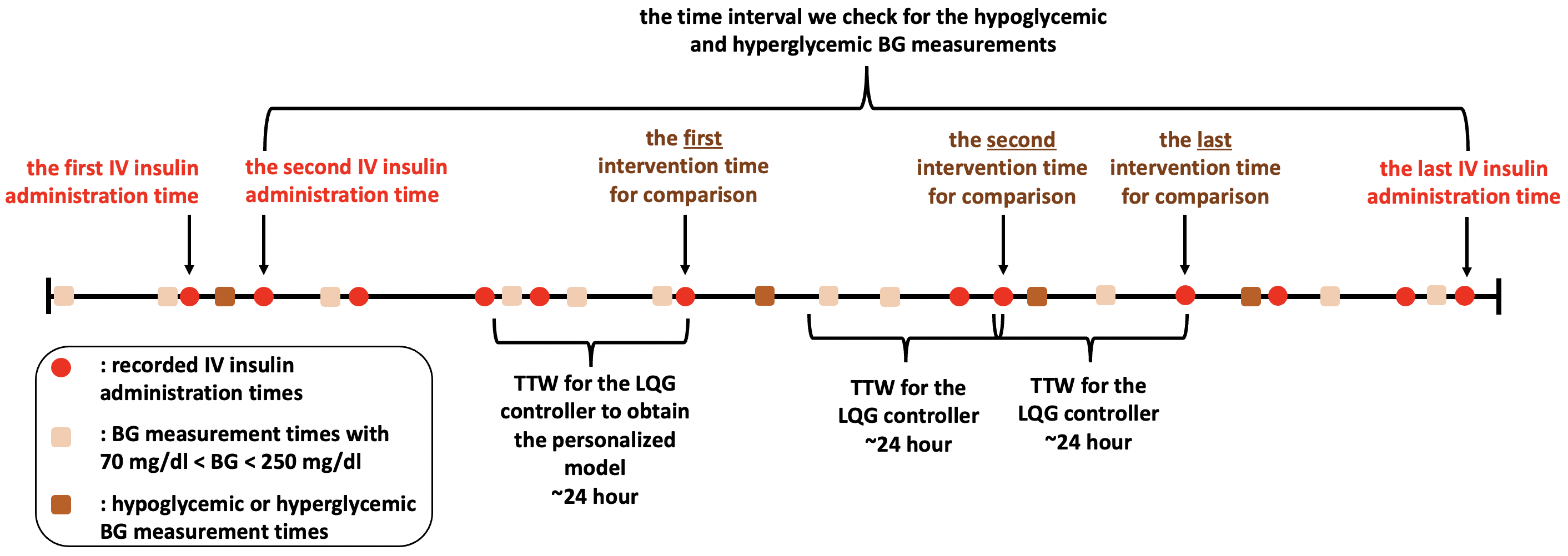}
	\caption{We show IV insulin administration and BG measurement times to describe our experimental design with the retrospective dataset. Recording times are shown in this figure only for illustration. Also, although this patient had four hypoglycemic/hyperglycemic measurements, we do not use the first one as it occurs before the second IV insulin administration time. So, there are three hypoglycemic/hyperglycemic measurements that we can use for evaluation of the GM approaches.}
	\label{real_data_schema}
\end{figure}

We need IV insulin to be administered over any TTW to quantify its effect on the patient through the MSG model. We need this information to estimate the optimal IV insulin rate using the LQG controller. So, we first define the time interval starting from the second to the last IV insulin administration times. We identify the times of hypoglycemic and hyperglycemic events over that time interval. Finally, the collection of the latest IV insulin administration times before each hypoglycemic or hyperglycemic measurement is the set of intervention times. The reason why we choose the intervention times before the hypoglycemic or hyperglycemic measurements is to evaluate and compare the efficiency of GM approaches in avoiding these hypoglycemic and hyperglycemic episodes. A visualization of this process is shown in Fig. \ref{real_data_schema} and an algorithmic description is given in Algorithm \ref{real_data_cont_alg}.

\subsection{Evaluation Techniques}\label{evaluation}

In this section, we describe the evaluation metrics that we use to compare the eﬀiciency of the LQG control and GM protocols in the simulated and real data settings. Note that the target region for the protocols are 120-150 mg/dL and 140-180 mg/dL. Since our goal is to evaluate the safety of the LQG controller and compare its eﬀiciency with these GM protocols, we adjusted the desired BG level for the controller for each comparison. We show the target region, hypoglycemic regions, and hyperglycemic regions for each protocol in Table \ref{regions_table}.

The first priority of GM in the ICU is to avoid life-threatening, severe hypoglycemia. However, whilst severe hyperglycemia is typically more manageable, GM approaches that avoid hyperglycemia are preferable to improve longer-term patient health status. Finally, maintaining BG levels in the target range is the most optimal outcome of GM as it is known to be the most beneficial for patients' health. We used the color scheme in Table \ref{regions_table} to represent the BG levels to be avoided or targeted. For example, a GM approach should avoid resulting BG levels in severe hypoglycemia and hyperglycemia (red) regions. The moderate (yellow) regions are not as risky but should still be avoided. While the target (green) region is the most preferable, the mild hypoglycemia and hyperglycemia (blue) regions are unlikely to pose any health risk for patients.

\begin{table}[!b]
	\centering
	\caption{Target region, hypoglycemic and hyperglycemic regions with the color codes used in the presentation of the experimental results. The BG values are in the unit of mg/dL.}
	\begin{tabular}{|l|c|c|}
		\hline
		\multicolumn{1}{|c|}{\textbf{Clinical Meaning}}				&	\textbf{Protocol A}				&	\textbf{Protocol B} \\ \hline
		\rowcolor{trred}
		Severe hyperglycemia	&	$(400, \infty)$										 & $(400, \infty)$ \\ \hline
		\rowcolor{tryellow}
		Moderate hyperglycemia	&	$(250,400]$										& $(250,400]$ \\ \hline
		\rowcolor{trblue}
		Mild hyperglycemia				&	$(150,250]$										& $(180,250]$ \\ \hline
		\rowcolor{trgreen}
		Target region							&	$[120,150]$										& $[140,180]$ \\ \hline
		\rowcolor{trblue}
		Mild hypoglycemia				&	$[70,120)$											&	$[70,140)$ \\ \hline
		\rowcolor{tryellow}
		Moderate hypoglycemia	&	$[40,70)$											& $[40,70)$ \\ \hline
		\rowcolor{trred}
		Severe hypoglycemia			&	$[0,40)$											& $[0,40)$	\\ \hline
	\end{tabular}
	\label{regions_table}
\end{table}

\subsubsection{Evaluation for Simulated Data}\label{eval_sim}

In this case, since we simulate the patient behavior by using ICUMM, we can ``deliver"  the LQG-suggested and protocol-suggested IV insulin rates to ``observe" their effect on the BG level. We denote the evaluation time points by $\{t_i\}_{i=1}^{P_k}$ for protocol and  $\{t_i\}_{i=1}^{L_k}$ for LQG controller where $P_k$ and $L_k$ are the number of evaluation times for the $k^{th}$ patient for $k = 1,2,...,N$. The BG measurements at these time points are denoted by $\{y^{protocol}_i\}_{i=1}^{P_k}$ and $\{y^{lqg}_i\}_{i=1}^{L_k}$, respectively, for $k = 1,2,...,N$. For each patient, we compute the minimum, maximum, and average of these BG measurements resulting from the IV insulin rates suggested by the protocol and LQG controller. Then we compare the resulting BG values by boxplots and hypothesis testing. Since we have resulting BG values for each patient by two different control approaches in this simulated data setting, we use paired sample t-test for comparing the mean difference in minimum, maximum, and average BG values. Note that the collection of minimum (resp. maximum) values provides the opportunity to compare the effectiveness of these two approaches for avoiding hypoglycemia (resp. hyperglycemia). Comparing average values helps us understand how well each approach keeps BG values in the target range on average.

\subsubsection{Evaluation for Real Data}\label{eval_real}

\begin{table}[t]
	\centering
	\caption{All the possible comparisons of $I^{real}$, $I^{protocol}$, and $I^{lqg}$ and their clinical meaning are shown.}
	\begin{tabular}{|p{8cm}|c|c|}
		\hline
		\multicolumn{1}{|c|}{\textbf{Clinical Outcome of Insulin Dosage}}		&	\textbf{Hypoglycemia}				&	\textbf{Hyperglycemia} \\ \hline
		\rowcolor{trblue}	&	$I^{lqg} = I^{protocol} < I^{real}$		&	\multicolumn{1}{|c|}{$I^{real} < I^{lqg} = I^{protocol}$}	\\
		\rowcolor{trblue}	&	$I^{lqg} < I^{protocol} < I^{real}$		&	\multicolumn{1}{|c|}{$I^{real} < I^{lqg} < I^{protocol}$}	\\
		\rowcolor{trblue}\multirow{-3}{*}{\parbox{8cm}{Both the LQG controller and the protocol give appropriate advice}}		&	$I^{protocol} < I^{lqg} < I^{real}$		&	\multicolumn{1}{|c|}{$I^{real} < I^{protocol} < I^{lqg}$}	\\ \hline
		\rowcolor{tryellow}		&	$I^{protocol} < I^{real} = I^{lqg}$		&	\multicolumn{1}{|c|}{$I^{lqg} = I^{real} < I^{protocol}$}	\\
		\rowcolor{tryellow}		& $I^{protocol} < I^{real} < I^{lqg}$		&	\multicolumn{1}{|c|}{$I^{lqg} < I^{real} < I^{protocol}$}		\\
		\rowcolor{tryellow}	\multirow{-3}{*}{\parbox{8cm}{The protocol gives more appropriate advice than the LQG controller}}		&	$I^{protocol} = I^{real} < I^{lqg}$		&	\multicolumn{1}{|c|}{$I^{lqg} < I^{real} = I^{protocol}$}		\\ \hline
		\rowcolor{trgreen}	&	$I^{lqg} < I^{real} = I^{protocol}$		&	\multicolumn{1}{|c|}{$I^{protocol} = I^{real} < I^{lqg}$}		\\
		\rowcolor{trgreen}	&	$I^{lqg} < I^{real} < I^{protocol}$		&	\multicolumn{1}{|c|}{$I^{protocol} < I^{real} < I^{lqg}$}		\\
		\rowcolor{trgreen}	\multirow{-3}{*}{\parbox{8cm}{The LQG controller gives more appropriate advice than the protocol}}	&	$I^{lqg} = I^{real} < I^{protocol}$		&	\multicolumn{1}{|c|}{$I^{protocol} < I^{real} = I^{lqg}$}		\\ \hline
		\rowcolor{trred}	&	$I^{real} < I^{lqg} = I^{protocol}$		&	\multicolumn{1}{|c|}{$I^{lqg} = I^{protocol} < I^{real}$}		\\
		\rowcolor{trred}	&	$I^{real} < I^{lqg} < I^{protocol}$		&	\multicolumn{1}{|c|}{$I^{lqg} < I^{protocol} < I^{real}$}		\\
		\rowcolor{trred}	\multirow{-3}{*}{\parbox{8cm}{Both the LQG controller and the protocol give inappropriate advice}}		&	$I^{real} < I^{protocol} < I^{lqg}$		&	\multicolumn{1}{|c|}{$I^{protocol} < I^{lqg} < I^{real}$}		\\ \hline
		\rowcolor{trpurple}	&	&	\multicolumn{1}{|c|}{}	\\
		\rowcolor{trpurple}		\multirow{-2}{*}{\parbox{8cm}{The LQG controller, protocol and real insulin rates are inappropriate}}
		&	\multirow{-2}{*}{$I^{lqg} = I^{protocol} = I^{real} > 0$}		&		\multicolumn{1}{|c|}{\multirow{-2}{*}{$I^{lqg} = I^{protocol} = I^{real}$}}		\\ \hline
		\rowcolor{trteal}	&	&	\multicolumn{1}{|c|}{}	\\
		\rowcolor{trteal}		\multirow{-2}{*}{\parbox{8cm}{The LQG controller, protocol and real insulin rates are appropriate}}
		&	\multirow{-2}{*}{$I^{lqg} = I^{protocol} = I^{real} = 0$}		&	\multicolumn{1}{|c|}{\multirow{-2}{*}{N/A}}	\\ \hline
	\end{tabular}
	\label{real_data_eval}
\end{table}

In this case, we use a different evaluation approach because of retrospective real-world data. We develop a comparison method for the protocol-suggested, $I^{protocol}$, and the LQG controller-suggested, $I^{lqg}$, rates through the IV insulin rate that was delivered to the patients, $I^{real}$. We perform this comparison only at the latest intervention times before moderate and severe hypoglycemic and hyperglycemic events. In this way, since we know the implication of $I^{real}$ (either hypoglycemia or hyperglycemia), by comparing $I^{protocol}$ and $I^{lqg}$ through $I^{real}$, we can test the efficacy of the LQG controller. Comparison of these three values could happen in different ways. which are grouped as shown in Table \ref{real_data_eval} for a meaningful comparison. Note also that the meaning of these comparison groups changes depending on the resulting BG measurement. For example, $I^{protocol} < I^{real} < I^{lqg}$ means that the LQG controller gives more appropriate advice than the protocol if the resulting BG is hyperglycemic, but the protocol gives more appropriate advice than the LQG controller if the resulting BG is hypoglycemic.

For the case when $I^{lqg} = I^{protocol} = I^{real}$, the value at which these IV insulin rates agree makes a difference for the clinical outcome if the resulting BG measurement is hypoglycemic. This is because if both the LQG controller and the protocol suggest no IV insulin administration and the real IV insulin rate is also 0 U/hr, even though the resulting BG is hypoglycemic, the advice is appropriate since they suggested the possible minimum IV insulin rate. On the other hand, if that value is larger than 0, this is considered inappropriate advice. Finally, since there is no definite upper bound for IV insulin rate, the case when $I^{lqg} = I^{protocol} = I^{real}$ is considered inappropriate regardless of their specific value when it results in hyperglycemic BG.

\section{Results}\label{num_res}

This section presents the results of comparing the LQG controller and protocol in the simulated and real data settings. Recall that we use two different GM protocols (protocol A and protocol B) in simulated and real-world data cases.

\subsection{Simulated Data Results}

We created virtual patients using ICUMM to evaluate and compare the effectiveness of the LQG controller to keep BG levels in the respective target range. We generated N = 200 patients and M = 20 different nutrition and IV insulin rate data. This means that there are 20 groups of patients, each consisting of 10 patients who received the same amount of nutrition over the entire time window (TTW and ETW) and IV insulin rate over their TTW. We set the lengths of TTW and ETW to be 24 hours each. We used the algorithmic structure (Section \ref{sim_data}) for each LQG controller-protocol pair and performed the comparison separately. The resulting plots of applying the LQG controller and one protocol are shown in the Supplementary Material.

\begin{subfigure}
	\setcounter{figure}{3}
	\setcounter{subfigure}{0}
	\centering
	\begin{minipage}[b]{0.5\textwidth}
		\includegraphics[width=\linewidth]{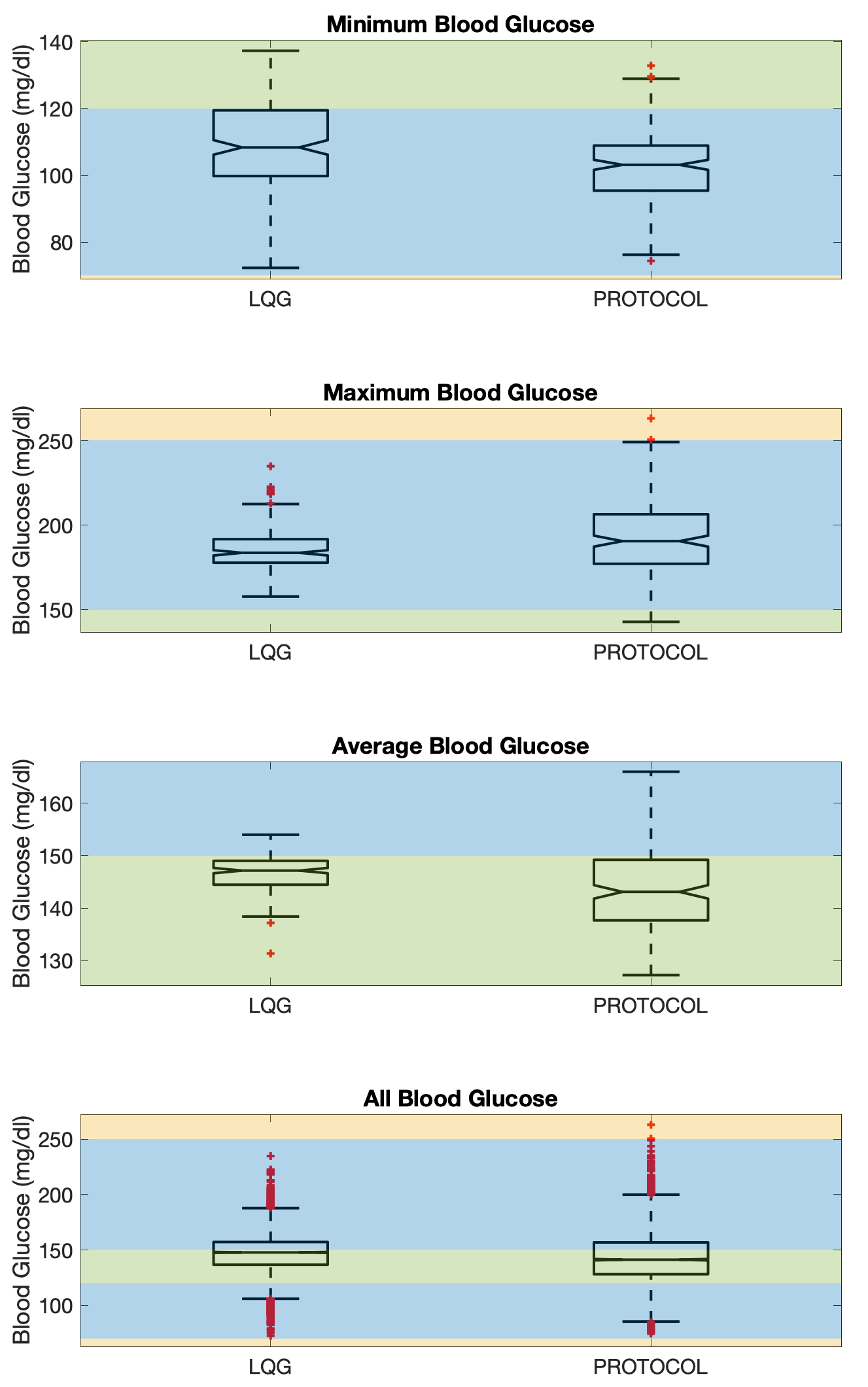}
		\caption{Protocol A}
		\label{fg:sub1}
	\end{minipage}%
	\setcounter{figure}{3}%
	\setcounter{subfigure}{1}%
	\begin{minipage}[b]{0.5\textwidth}
		\includegraphics[width=\linewidth]{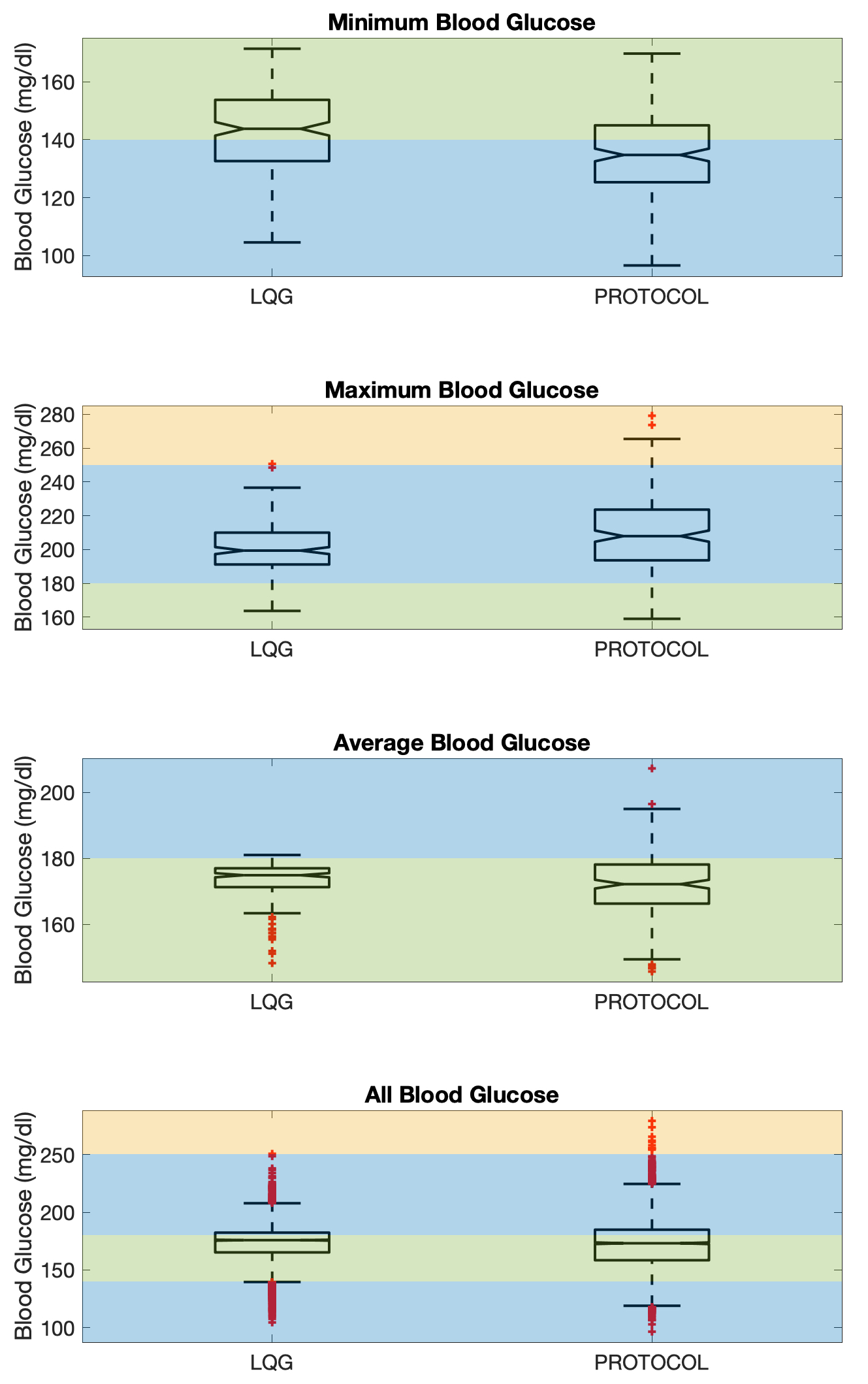}
		\caption{Protocol B}
		\label{fig:sub2}
	\end{minipage}%
	\setcounter{figure}{3}
	\setcounter{subfigure}{-1}	
	\caption{The figure shows minimum, maximum, average, and all BG values ``measured" resulting from an intervention for all virtual patients. Minimum (resp. maximum) BG shows if the patient had any hypoglycemic (resp. hyperglycemic) event. The average BG shows the overall performance of a GM approach in keeping BG levels in the target range. All BG shows all the BG values ``measured" for all patients. LQG controller is tuned for the respective protocol's target region, using the same intervention times with each protocol together with the nutrition change times. This color scheme is described in Section \ref{evaluation}.}
	\label{sim_data_boxplot}
\end{subfigure}

We performed hypothesis testing to quantify the statistical significance of the difference between the resulting minimum, maximum, and average BG levels for each patient by the pairs of LQG controller and the protocols. Our null hypothesis is that there is no difference between the minimum (resp. maximum, average) BG values. The alternative hypothesis is that the difference between the minimum (resp. maximum, average) BG values is nonzero. We used MATLAB's \texttt{ttest} function to test these hypotheses. We checked the normality assumption by using MATLAB's \texttt{kstest} function for all samples. The test failed to reject the null hypothesis stating that the sample comes from a standard normal distribution, against the alternative that it does not come from such a distribution, at $\alpha = 0.05$ significance level for all three samples of minimum, maximum, and average BG values (respective p-values are $0.27, 0.61, 0.34$). Then, using \texttt{ttest}, we concluded the statistical significance of the differences between the minimum (resp. maximum, average) BG values under the significance level $\alpha = 0.05$ (respective p-values are $1.2*10^{-9}, 5.7*10^{-11}, 0.03$). We show the 95\% CIs for these paired differences in Table \ref{sim_data_hyp_testing}.

Our results show that the LQG controller provided more appropriate IV insulin dose recommendations than the GM protocols resulting in a great fraction of measured BG levels in, or close to, the target region. All three GM approaches avoided severe and moderate hypoglycemia and severe hyperglycemia for all virtual patients. Moreover, the LQG controller avoided moderate hyperglycemia for all virtual patients.

\textit{First}, even though all three GM approaches avoided severe and moderate hypoglycemia (below 70 mg/dL) at any time over the ETW, the LQG protocol results in slightly higher minimum BG levels that are closer to the target region (first-row panels of Fig. \ref{sim_data_boxplot}). The 95\% CIs for the paired difference of minimum of the resulting BG values (the first row of Table \ref{sim_data_hyp_testing}) shows that LQG controller maintains BG values further from the hypoglycemic region than both protocols.

\textit{Second}, while all three GM approaches can avoid severe hyperglycemia, the LQG controller also avoids moderate hyperglycemia as shown in the second-row panels of Fig. \ref{sim_data_boxplot}. In the LQG controller - protocol B comparison (the second-row right-hand panel of Fig. \ref{sim_data_boxplot}), we see that the LQG controller results in only one maximum BG value (250.7 mg/dL) on the boundary between mild and moderate hyperglycemia and all the remaining maximum BG values lie either in the mild hyperglycemia region or the target region. In addition, the 95\% CIs for the difference between the paired maximum BG values (the second row of Table \ref{sim_data_hyp_testing}), when considered together with the resulting BG values (the second-row panels of Fig. \ref{sim_data_boxplot}), show that the difference in avoiding severe and moderate hypoglycemia is statistically significant, i.e., the LQG controller results in BG values further from the severe and moderate hyperglycemia regions than both protocols.

\textit{Third}, the LQG controller maintains a larger number of average BG values in the target region compared to protocols (third-row panels of Fig. \ref{sim_data_boxplot}). The average BG levels over the ETW for each patient obtained by an LQG controller-suggested IV insulin rate mostly lie on the respective target region for each protocol. The 95\% CIs for the paired difference of average BG values (third row of Table \ref{sim_data_hyp_testing}) show that the difference is statistically significant. Notably, the LQG controller, when using protocol B-suggested intervention times, resulted in almost all BG levels in the target region. This supports that the LQG controller indeed learns the patients' BG dynamics through the relatively simple yet still physiological BG model, the MSG model, and suggests appropriate IV insulin rates.

\begin{table}[!ht]
	\centering
	\caption{This table shows the 95\% CIs for the paired sample differences of the minimum, maximum, and average BG values resulting from LQG controller and protocol, computed separately for each virtual patient over their ETW. The LQG controller provides more appropriate IV insulin administration recommendations resulting in BG levels in or closer to the respective target ranges (CIs for minimum and maximum BG).}
	\begin{tabular}{|c|c|c|}
		\hline
		&	\textbf{LQG \& Protocol A}		&	\textbf{LQG \& Protocol B}	\\ \cline{2-3}
		& \multicolumn{2}{c|}{$\text{BG}^{lqg}-\text{BG}^{protocol}$}	 \\ \hline
		Minimum Blood Glucose		&	[3.5,8.4]			&	[5.6,10.7]	\\ \hline
		Maximum Blood Glucose		&	[-9.4,-4.2]			&	[-10.0,-5.6]	\\ \hline
		Average Blood Glucose			&	[1.7,4.5]			&	[0.2,2.4]	\\ \hline
	\end{tabular}
	\label{sim_data_hyp_testing}
\end{table}

\textit{Finally}, the fourth-row panels of Fig. \ref{sim_data_boxplot} show the collection of all BG values over the ETW for all 200 patients. These boxplots show that the collection of BG values obtained as the result of the LQG controller-suggested IV insulin rates cover a slightly narrower BG range compared to each protocol. In this case, since the LQG controller is developed to account for the nutrition change times, the resulting BG values cannot be paired and we cannot perform similar hypothesis testing to quantify the statistical significance of the difference.

These results show that the mechanistic model-based personalized GM approach provides efficient IV insulin administration to avoid adverse events and maintains a larger fraction of BG values in or closer to the desired target regions.

\begin{table}[t]
	\centering
	\caption{We show the performance of the pairs of LQG controller - protocol using retrospective real-world data. We use the same dataset for both LQG controller - protocol pair. The LQG controller gave more appropriate advice than protocol A to avoid hypoglycemia in most of the cases. When compared with protocol B, it is not possible to distinguish the effectiveness. For the hyperglycemic events, the LQG controller gave more appropriate advice than both protocols in most of the cases.}
	\begin{tabular}{|p{6cm}|c|c|c|c|}
		\hline
		\multicolumn{1}{|c|}{\multirow{2}{*}{\parbox{6cm}{\textbf{Clinical Outcome of Insulin Dosage}}}}				 &	\multicolumn{2}{|c|}{\textbf{Hypoglycemia (19 events)}}				&	\multicolumn{2}{|c|}{\textbf{Hyperglycemia (107 events)}}\\ \cline{2-5}
		&		\textbf{Protocol A} 		&\textbf{Protocol B}		&	\textbf{Protocol A} &\textbf{Protocol B}	\\ \hline
		\rowcolor{trblue}	&	&	&	&	\\
		\rowcolor{trblue}	\multirow{-2}{*}{\parbox{6cm}{Both the LQG controller and the protocol gave appropriate advice}}	&	\multirow{-2}{*}{6}	&	\multirow{-2}{*}{16}	&	\multirow{-2}{*}{23}	&	\multirow{-2}{*}{20} \\ \hline
		\rowcolor{tryellow}	&	&	&	&	\\
		\rowcolor{tryellow}	\multirow{-2}{*}{\parbox{6cm}{The protocol gave more appropriate advice than the LQG controller}}	&	\multirow{-2}{*}{0}	&	\multirow{-2}{*}{0}	&	\multirow{-2}{*}{18}	&	\multirow{-2}{*}{19} \\ \hline
		\rowcolor{trgreen}	&	&	&	&	\\
		\rowcolor{trgreen}	\multirow{-2}{*}{\parbox{6cm}{The LQG controller gave more appropriate advice than the protocol}}	&	\multirow{-2}{*}{13}	&	\multirow{-2}{*}{0}	&	\multirow{-2}{*}{56}	&	\multirow{-2}{*}{52} \\ \hline
		\rowcolor{trred}	&	&	&	&	\\
		\rowcolor{trred}	\multirow{-2}{*}{\parbox{6cm}{Both the LQG controller and the protocol gave inappropriate advice}}	&	\multirow{-2}{*}{0}	&	\multirow{-2}{*}{0}	&	\multirow{-2}{*}{8}	&	\multirow{-2}{*}{12} \\ \hline
		\rowcolor{trpurple}	&	&	&	&	\\
		\rowcolor{trpurple}	\multirow{-2}{*}{\parbox{6cm}{The LQG controller, protocol and real insulin rates were inappropriate}}	&	\multirow{-2}{*}{0}	&	\multirow{-2}{*}{0}	&	\multirow{-2}{*}{2}	&	\multirow{-2}{*}{4} \\ \hline
		\rowcolor{trteal}	&	&	&	&	\\
		\rowcolor{trteal}\multirow{-2}{*}{\parbox{6cm}{The LQG controller, protocol and real insulin rates were appropriate}}	&	\multirow{-2}{*}{0}	&	\multirow{-2}{*}{3}	&	\multirow{-2}{*}{N/A}	&	\multirow{-2}{*}{N/A} \\  \hline
	\end{tabular}
	\label{real_data_tab}
\end{table}

\subsection{Real Data Results}

Now, we use retrospective data to evaluate the performance of the LQG controller. Unlike a prospective study, with retrospective data, we cannot observe the direct implications of alternate GM strategies on patients. This limitation shaped our evaluation approach as discussed in Section \ref{eval_real}.  We compared the eﬀicacy of the LQG controller and the protocols based on their suggestions at the latest intervention times before the recorded hypoglycemic and hyperglycemic events.

Our results show that the LQG controller recommended IV insulin administration rates at least as effective as the protocol in avoiding hypoglycemic and hyperglycemic events. This result is reflected in the values shown in Table \ref{real_data_tab}, which were obtained by the comparison of the LQG controller-suggested ($I^{lqg}$) and protocol-suggested ($I^{protocol}$) IV insulin rates through the real IV insulin rate ($I^{real}$) as shown in Table \ref{real_data_eval}. Since we compared the LQG controller with each protocol separately and the relative difference between IV insulin values affects the resulting values in Table \ref{real_data_tab}, we only compare the performance of the LQG controller with the protocols separately rather than comparing the eﬀicacy of the protocols.

\textit{First}, the LQG controller provided more appropriate advice than protocol A for 6 out of the 19 hypoglycemic measurements. This shows that with the IV insulin rate suggested by the LQG controller, a large portion of the hypoglycemic events could be potentially avoided. On the other hand, when compared with protocol B, 16 out of 19 times both the LQG controller and the protocol gave useful advice to avoid the larger portion of the hypoglycemic events. For these cases, it is impossible to know which (LQG controller or protocol) was more accurate because the data are retrospective. We also note that for 3 out of 19 hypoglycemic events, the three IV insulin rates ($I^{lqg}$, $I^{protocol}$, and $I^{real}$) agreed at the rate of 0, which still caused hypoglycemia. So, we conclude that the suggested rates were still reasonable since it is not possible to suggest an IV insulin rate that could avoid hypoglycemia.

\textit{Second}, the LQG controller gave more helpful advice than protocol A at 56 and protocol B at 52 out of all 107 hyperglycemic measurements. This is because we set the desired BG level to be the upper bound of the respective target region for each protocol to avoid hypoglycemic episodes. When we consider the distribution of the values of all five possible outcomes for hyperglycemic events, we see that the event LQG controller providing more appropriate IV insulin rate suggestion than both protocols is the most likely outcome. Around 18-19 out of 107 hyperglycemic events, the protocols outperformed the LQG controller, and 20-23 out of the total 107 hyperglycemic events, both the protocols and the LQG controller gave reasonable advice that could potentially avoid these events.

\textit{Third}, the value of events that both the LQG controller and protocol gave inappropriate advice is around 8-12 times out of all the hyperglyceic events, which indicates that all the GM approaches could be further improved to avoid these adverse events.

These results show that a model-based GM strategy can potentially avoid all hypoglycemic events and reduce the incidence of hyperglycemic events to much lower rates by providing personalized IV insulin administration suggestions.

\textit{In summary}, these results exhibit the potential improvement that the LQG controller we developed based on the MSG model could improve the GM outcomes by avoiding a large portion of hypoglycemic and hyperglycemic events.

\section{Discussion}\label{discussion}

\subsection{Overview of Findings} We developed a model-based GC approach accounting for the limitations of routinely collected ICU data. The controller requires (1) a mechanistic model describing BG behavior and (2) a control algorithm. Because of the aforementioned challenges, such as unmodeled interventions and sparsely collected data, ICU patients' BG behavior cannot be resolved accurately enough to be used for model-based control methodologies. To address these challenges, we used the MSG model, a linear and stochastic BG model. The model's simplicity allowed us to overcome data sparsity limitations and to track the mean BG behavior accurately while its stochasticity helped quantify phenomena excluded from the model and, in particular, resolve the BG variance providing useful insight for GM strategies. Then, we developed an LQG controller based on the MSG model. We designed simulated and real-world retrospective data experiments to evaluate the performance of this methodology by comparing its effectiveness with operational GM flow charts. By mimicking the conditions of a prospective study, the simulated data experiment provided a ground for a fair comparison of different GM approaches. The real-world retrospective data experiment provided an opportunity to understand how well the GM approaches can handle real-world data limitations. Our results showed that it is possible to implement an efficient, actionable, and safe model-based GM technique based on a simple stochastic mechanistic BG model paired with a suitable control algorithm.

\subsection{Potential Benefits of Using the MSG Model} We used the MSG model as the basis of this control theoretical GM approach to exploit the computational advantages of such a simple and analytically solvable yet physiological model. A linear (stochastic) model with an analytical solution is advantageous because (1) it reduces computation time significantly, making this control framework applicable in an online setting as a CDS tool; (2) it avoids errors affecting predictive performance due to using numerical approximations when an analytic solution does not exist; and (3) it allows us to use a simple and eﬀicient control technique, LQG controller, enhancing the computational effectiveness of the whole framework.

\subsection{Differences Between the Protocols and Their Effect on the LQG Controller}  Our results with simulated data showed that both protocols resulted in similar GM performance. One significant difference between the protocols is their performance in avoiding lower BG values. Even though both protocols place BG values in safe regions, the \textit{Minimum Blood Glucose} plots of Fig. \ref{sim_data_boxplot} show that protocol B results in a larger number of minimum BG values in or closer to its target region. Since the LQG controller uses the corresponding intervention times for each protocol (in addition to patients' nutrition change times), the difference between the suggested intervention times by the protocols might be a contributing factor in avoiding lower BG values. We will investigate the reasons causing this difference in the development of a more effective GM controller.

\subsection{Different Insulin Types and Delivery Routes in the ICU} Insulin for GM in the ICU is delivered intravenously during acute illness. As their health status improve, they are transitioned to subcutaneous insulin. IV insulin is typically only short-acting. However, depending on the needs of patients, subcutaneous insulin could be rapid-acting, short-acting, intermediate-acting, long-acting, or ultralong-acting. We developed this methodology estimating the optimal IV insulin amount to evaluate the applicability of in the ICU setting. Including the subcutaneous insulin with all possible different types in this methodology will pose additional complexity, which will be addressed in a future work.

\subsection{Results of Our Conservative Approach in the Controller Development} The 95\% CIs for average and all BG values for virtual patients, shown in Table \ref{sim_data_hyp_testing}, show that all the CIs for the LQG controller represent a higher range than the corresponding CI of the protocol. This behavior results from our conservative approach, i.e., setting the target value to be the upper bound of the respective target range and bounding IV insulin rate by 25 U/hr from above for the LQG controller. This result indicates that our specific constraints in developing the LQG controller are reflected in the GM results.

\subsection{Potential Benefits of Personalized Glycemic Management} A common behavior observed in Fig. \ref{sim_data_boxplot} is that the LQG controller results in less variability than protocols in the BG values. Unlike the protocols, which follow a one-size-fits-all GM approach, the LQG controller provides personalized interventions based on the information extracted from patients' data via the MSG model. Personalized GM toward the same target value for all patients might be the reason for reduced variability in the resulting BG values. Importantly, a GM approach could be even more personalized by setting a unique target range for each patient. Such a controller is beyond the scope of the approach we develop here and requires investigation by clinicians. However, a patient-specific BG range is very straightforward to include in this control algorithm. This easily adjustable feature of the model-based control approaches makes them practical and desirable within CDS tools once they are externally validated.

\subsection{Potential Benefits of Virtual Patient Cohorts for Evaluation} Creating a cohort of virtual patients for the initial investigation of an intervention or treatment strategy provides a useful platform for evaluating the safety of the intervention and for a fair comparison of different intervention strategies. We used of these features in the evaluation of our control approach. For example, while the amount of IV insulin delivered to patients in the real ICU setting is constrained by situational factors, with virtual patients we can simulate the delivery of any amount suggested by the LQG-controller or a protocol without worrying about safety. This showed us the possible range of resulting BG levels. Also, comparing two different intervention strategies on virtual patients mimics the case of comparing two intervention strategies on two identical patients, which is impossible in any real-world setting. In this way, we aim to isolate the effect of the interventions from other factors.

\subsection{Generalizability} The underlying mathematical approach of this GM methodology shows some similarities with the AP project. However, there are some fundamental differences between this methodology and the AP. The AP is developed to deliver insulin automatically via an insulin pump based on continuous glucose monitoring (CGM) data for GM of T1DM patients. This methodology can work with CGM data, but does not require that data. Also, the AP does not typically incorporate the nutritional intake of patients into the model-based control system. This methodology was not developed for automatic insulin delivery but for use as a CDS tool for GM of ICU patients regardless of their diabetic state, which is captured and accounted by patient-specific model estimation. In addition, the way this methodology accounts for ICU-specific challenges, such as the rapidly changing BG behavior and data sparsity can be transferred to other ICU settings with similar limitations. Also, the control pipeline is not specific to GM. It could easily be transferred to other healthcare settings requiring clinical decision support, such as optimal timing of drug delivery for oncology patients \citep{martinez2021clinical}.

\subsection{Limitations and Future Work} The methodology is developed based on the MSG model. While this model has advantages for GC purposes in the ICU setting, it also has some limitations. The mean component of the stochastic model output only provides information about the mean behavior of the system and generally cannot forecast the complex dynamical behavior. In this sense, the MSG model might oversimplify the system. We loose the ability to track the exact trajectory by using this stochastic MSG model. Therefore, if the exact details of the complex system behavior are important to capture, we cannot do this by the MSG model.

Using retrospective data poses certain limitations to the effective evaluation of GM strategies. Unlike a prospective study, it is impossible to observe the effect of suggested interventions. However, it still provides useful insight into the effectiveness of the GC approaches. Another limitation of the real-world retrospective dataset is the total number of patients included in it. Note that the dataset contains all patients satisfying the inclusion criteria described in Section \ref{data}. This limitation is because of the current version of the methodology accounting for IV insulin delivery and tube fed nutrition as it restricts the number of patients that can be included in the dataset. Additionally, our evaluation method greatly reduced the number of patients because it required the presence of an adverse event. For example, most patients require subcutaneous insulin in addition to or instead of IV insulin. These patients were not included in our study. Including these cases is a goal of future work. In addition, the physiological variants that occur 10\% (resp. 5\%) of the time have a 91\% (resp. 69\%) chance of being represented in this subgroup with 23 patients. Therefore, while not capturing all possible physiological variants, we would get most of the common human variants. Also, using retrospective data for evaluation prevents measuring improvements in patient outcomes, which can be possible with prospective studies. The current version of this methodology is not ready for prospective studies since it needs to be expanded to be used for a broader ICU patient population, which will be addressed in future studies.

We used ICUMM to generate the profile and data of virtual patients, as described in Section \ref{sim_data}. We generated model parameters independently using a uniform distribution over the given range for each parameter to represent the patient profiles. A more realistic scenario would be to account for the correlation between these parameters. Unfortunately, the correlation level, in the form of a covariance matrix, to generate more realistic virtual patient profiles is unknown. Further, setting the model parameters constant while generating time-varying ICUMM parameters could better reflect the non-stationary BG behavior of ICU patients. We will address these limitations in future work.

Other future studies will focus on expansion and improvement of this methodology to the broader ICU patient population. Exogenous insulin is delivered through different routes (IV and subcutaneous) and in different types (e.g., rapid-acting, short-acting, and long-acting). We will incorporate subcutaneous delivery component and the impact of other insulin types into the model. Also, we will adjust the model to account for the effect of orally consumed nutrition on BG levels. The LQG control algorithm will immediately be applicable to estimate optimal exogenous insulin amount after these modifications in the model. Therefore the methodology will be available to provide GM decision support for most of the ICU patients. Also, some medications used in ICU (e.g., glucocorticoids, fluoroquinolones, and beta-blockers)  are known the affect BG dynamics of patients, but these effects are not accounted for in the GM protocols, as the protocols suggest insulin dose adjustments after observing the changes in the BG levels. We will incorporate the effect of these medications into the MSG model to be adjust exogenous insulin delivery in advance to avoid BG levels deviate from the target BG region. The LQG controller requires at least one insulin administration data for training, which might be considered a limitation. However, with any model-based controller, one needs to ``learn" the effect of the exogenous insulin on BG dynamics through the estimated model parameters and reflect this information to the predicted amount of optimal exogenous insulin. In practice, the LQG controller would be activated once the first insulin dose is given according to a local protocol. Finally, we will use qualitative methodology tools to investigate how we can integrate these quantitative GM techniques into the clinical workflow for efficient and actionable CDS tools. Completion of these modeling studies will position us to design prospective studies to evaluate the effectiveness of this methodology for GM in the ICU and for improving patient outcomes. For example, we will perform a multi-site clinical trial with healthcare professionals on the control group following the state-of-the-art GM protocol flowchart and the healthcare professionals on the experimental group being aided by the clinical decision support developed based on the resulting methodology.

\subsection{Conclusion}\label{conc}

We investigated the effectiveness of a model-based controller developed using a physiology-based linear stochastic model representing BG dynamics in the ICU setting. The BG behavior of ICU patients are highly non-stationary due to these patients' critical illness, constant tube feeding, and frequent interventions. In addition, the routinely collected data in ICU is sparse, making it more difficult to use for control purposes. This model-based GM approach provides real-time safe and efficient intervention strategies using routinely collected data. Moreover, the methodology is easily adjustable to achieve other BG target ranges for different GM strategies. These features of this methodology make it usable for CDS tools within EHR systems. Also, the general structure of this control framework is transferable to other ICU settings, such as ventilation and hemodynamic management that pose similar characteristic challenges. We will address the limitations listed in Section \ref{discussion} before an extensive evaluation in a pilot study.

\section*{Conflict of Interest Statement}

The authors declare that the research was conducted in the absence of any commercial or financial relationships that could be construed as a potential conflict of interest.

\section*{Author Contributions}

M.S., A.M.S., and D.J.A. conceptualized this work. Y.W. and J.N.S. curated the data. G.H., A.M.S., and D.J.A. acquired funding to perform this work. M.S. designed the methodology, performed the formal analysis, developed the software, created figures and tables, and performed validation with supervision from G.H., C.C.L.W., T.D.B., A.M.S., and D.J.A. Also, M.S. wrote the first draft. All the authors read, reviewed, and edited the manuscript.

\section*{Funding}

This work was supported in part by NIH R01 LM006910 ``Discovering and Applying Knowledge in Clinical Databases" and R01 LM012734 ``Mechanistic Machine Learning."

\section*{Acknowledgments}

A.M.S. is grateful for support from the Department of Defense Vannevar Bush Faculty Fellowship.

\section*{Data Availability Statement}

The current version of the data that we used here are fully identifiable, and we will not share this version to protect patient and healthcare professional privacy. We are currently working to create a subset of non-PHI, de-identified, hand-verified data publicly available (subject to agreements with the institutions where the data originate) on data sharing repositories.

\bibliographystyle{Frontiers-Harvard} 
\bibliography{main}


\end{document}


\onecolumn
\firstpage{1}

\title {{\helveticaitalic{Supplementary Material}}}

\maketitle

\section{Data}

Our retrospective evaluation methodology required including patients who experienced at least one adverse event (either hypoglycemia or hyperglycemia), which resulted in 23 patients. However, the original dataset included 106 patients who required IV insulin administration and fed through an enteral tube. It is important to note that 20 out of the 106 patients required only bolus insulin which was delivered intravenously. The remaining 86 patients required either only continuous IV insulin or bolus IV insulin in addition to continuous IV insulin. Here, we provide detailed information about all of the 106 patients in Table \ref{table1}. Since the current version of the methodology provides decision support for continuous IV insulin delivery, we provide the length of continuous IV insulin delivery statistics only for the 86 patients.

\begin{table}[!ht]
	\centering
	\caption{Demographic and ICU stay-related health record information of patients who required IV insulin administration and fed through an enteral tube in the real-world dataset. These patients did not necessarily experience an adverse event (hypoglycemia or hyperglycemia). (* indicates that the reported statistics represent 86 patients who required continuous IV insulin administration. The remaining 20 patient were delivered only bolus IV insulin.)}
	\begin{tabular}{| p{0.45\linewidth} | r| r | r | r | r |}
		\hline
		& \textbf{mean} & \textbf{stdev} & \textbf{min} & \textbf{median} & \textbf{max}\\ \hline
		\textbf{Age (years)}	& 55.6 		&	13.9				&	25.0			&	56.0					& 81.0\\ \hline
		\textbf{Length of ICU stays that contain an interval meeting the inclusion criteria (days)}	& 16.3	&	15.2	&	3.4	&	10.8	&	72.9\\ \hline
		\textbf{Total number of BG measurements}	&	89.1	&	75.3	&	4	&	72.5	&	462\\ \hline
		\textbf{Length of tube feed administration (hours)}	& 194.7	&	154.9	&	73.3	&	142.8	&	875\\ \hline
		\textbf{*Length of IV insulin administration (hours)}	& 96.7	&	82.8	&	1.5	&	84.5	&	444.0\\ \hline
		\textbf{Race}		& \multicolumn{5}{l|}{American Indian and Alaska Native: 2}\\
		& \multicolumn{5}{l|}{Asian: 4}\\
		& \multicolumn{5}{l|}{Black or African American: 12}\\
		& \multicolumn{5}{l|}{Multiple Race: 1}\\
		& \multicolumn{5}{l|}{Other: 6}\\
		& \multicolumn{5}{l|}{White or Caucasian: 81}\\ \hline
		\textbf{Ethnicity}		& \multicolumn{5}{l|}{Hispanic: 12}\\
		& \multicolumn{5}{l|}{Non-Hispanic: 94}\\ \hline
		\textbf{Sex}	& \multicolumn{5}{l|}{Female: 36}\\
		& \multicolumn{5}{l|}{Male: 70}\\ \hline
	\end{tabular}
	\label{table1}
\end{table}

\section{Simulations of Blood Glucose Levels Resulting from the Model-Based Controller and Protocol Used in Clinical Practice}

We show the simulations of resulting BG values with the LQG controller-suggested and protocol-suggested IV insulin rates for a detailed explanation of the developed methodology. In Fig. \ref{simulation}, the top panels show simulated BG values, representing a virtual patient's response to the IV insulin rate shown in the middle panels. The nutrition rates (lower panels) are the same in each subfigure.

We used the ICUMM to create virtual patient profiles and to represent those virtual patients' glycemic response to LQG-suggested and protocol-suggested IV insulin rates. In this sense, ICUMM is a representation of reality. We used the MSG model in our control-theoretical framework to \textit{learn} the patients' BG dynamics and estimate the optimal amount of IV insulin rate for GM. This figure demonstrates potential GM performance given a perfect representation of BG dynamics with a mechanistic model. The target value for the LQG controller is the upper limit of the respective protocol's target range (180 mg/dL) to avoid hypoglycemia. The BG values simulated with the MSG model hover around that target value, with occasional lower values due to nutrition rate changes shifting the system equilibrium. However, the LQG controller quickly identifies the IV insulin rate to shift the equilibrium to the target value of 180 mg/dL.

\begin{subfigure}[ht]
	\setcounter{figure}{1}
	\setcounter{subfigure}{0}
	\centering
	\begin{minipage}{.33\textwidth}
		\centering
		\includegraphics[width=\linewidth]{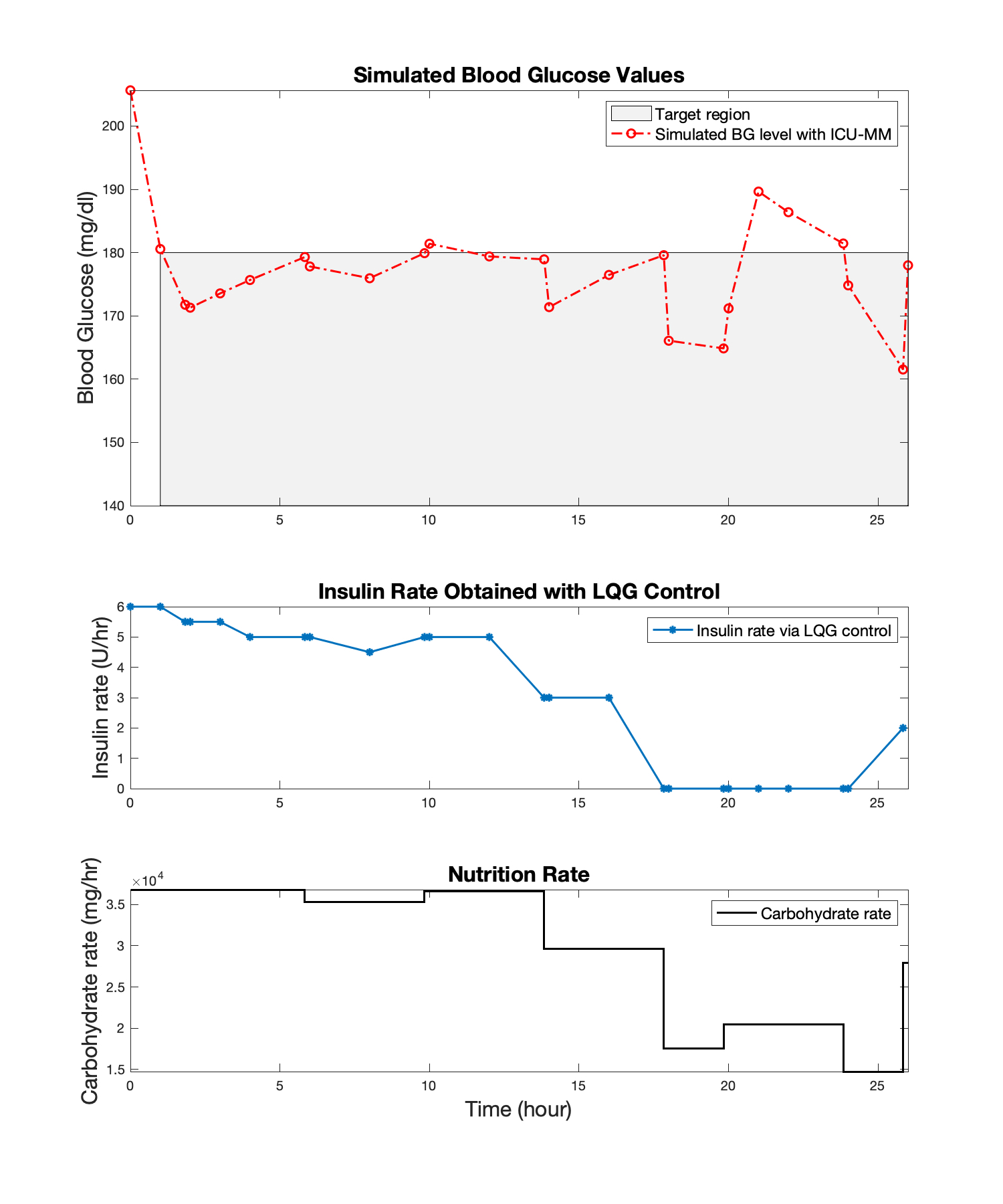}
		\caption{LQG --- ICUMM}
		\label{lqg_icumm_sim}
	\end{minipage}%
		\setcounter{figure}{1}%
	\setcounter{subfigure}{1}%
	\begin{minipage}{.33\textwidth}
		\centering
		\includegraphics[width=\linewidth]{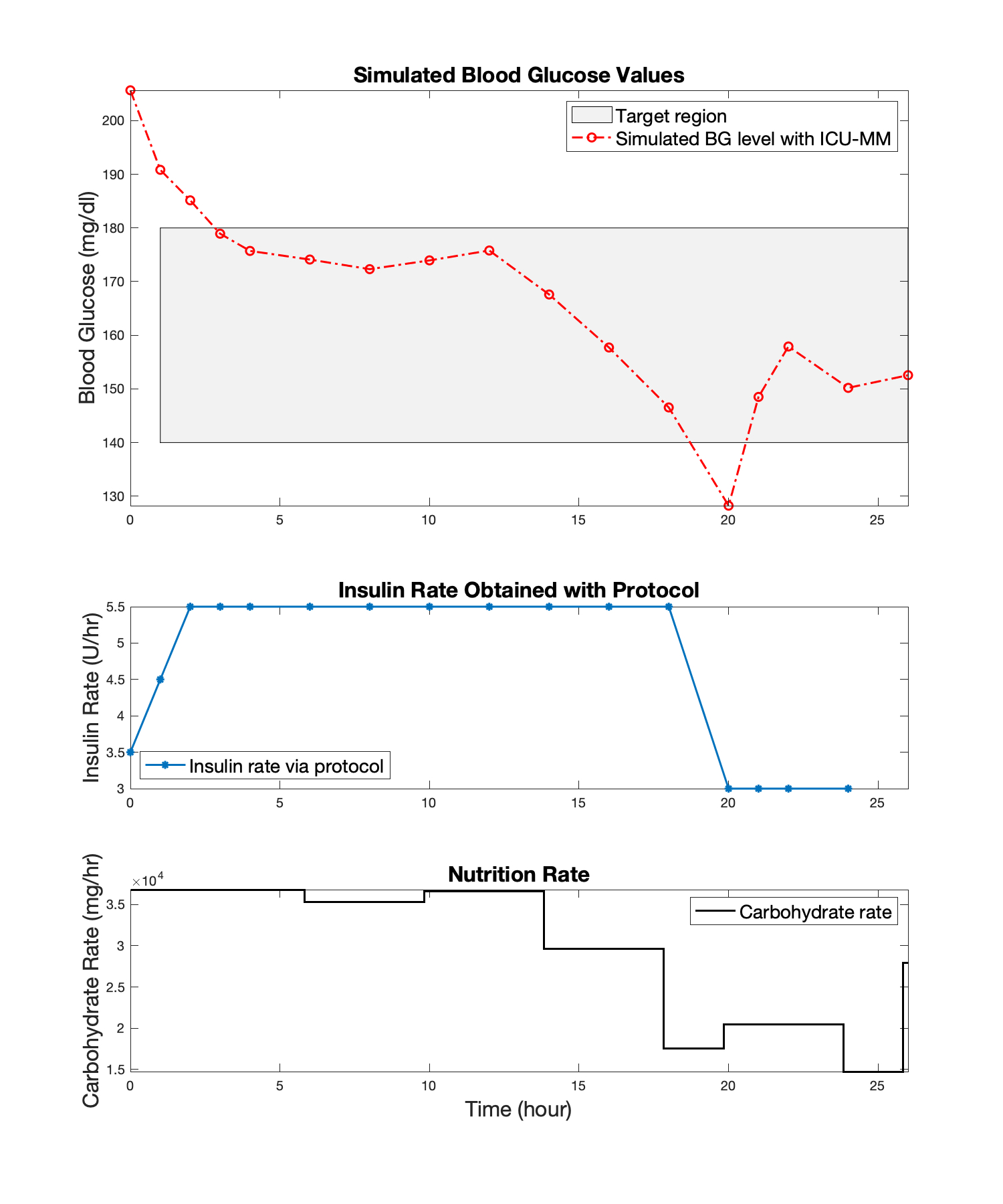}
		\caption{Protocol --- ICUMM}
		\label{protocol_sim}
	\end{minipage}
		\setcounter{figure}{1}%
	\setcounter{subfigure}{2}%
	\begin{minipage}{.33\textwidth}
		\centering
		\includegraphics[width=\linewidth]{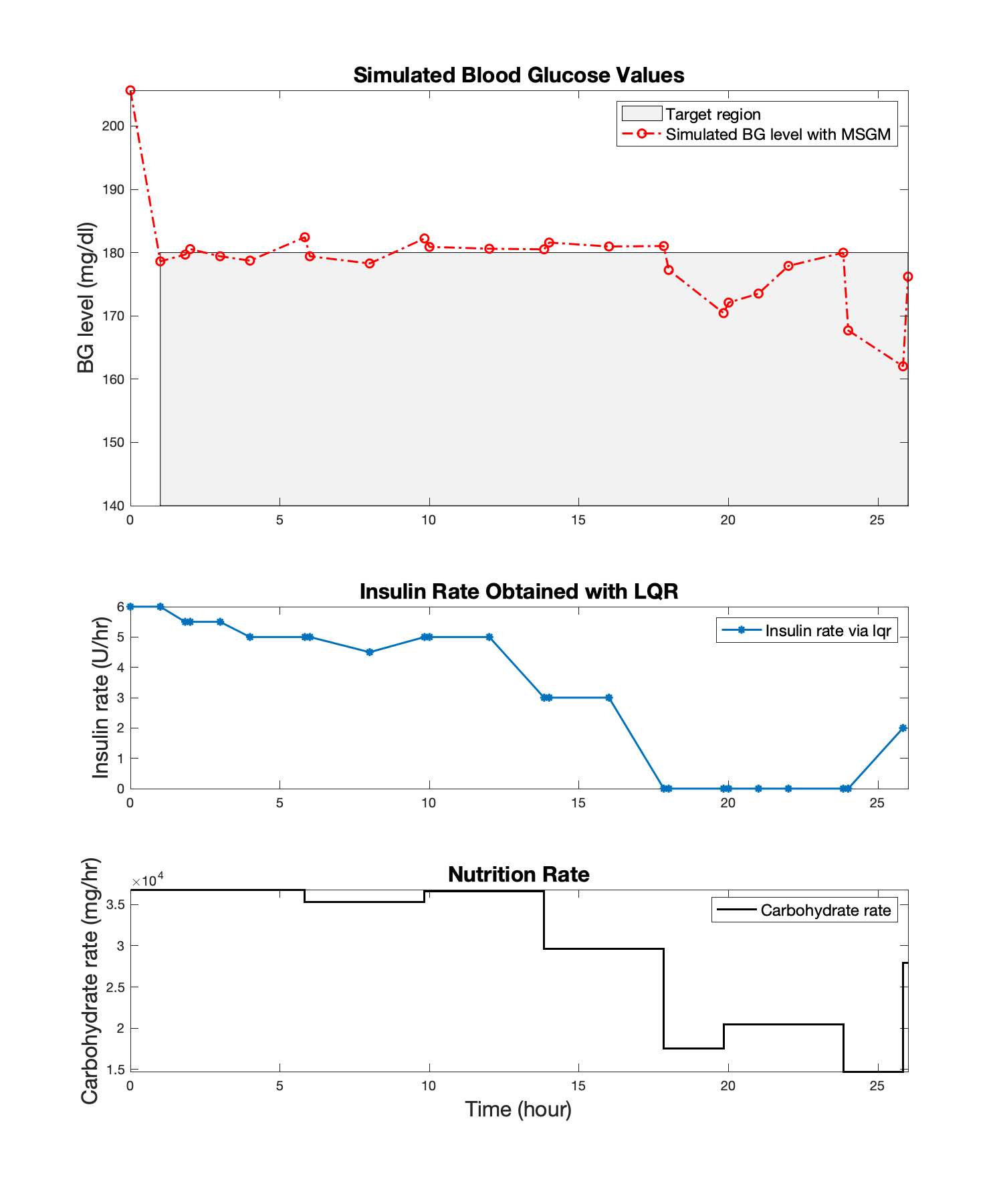}
		\caption{LQG - MSGM}
		\label{lqg_msgm_sim}
	\end{minipage}
	\setcounter{figure}{1}
	\setcounter{subfigure}{-1}	
	\caption{Tthe upper panels show the BG values simulated by the model shown in the respective legend in response to the IV insulin rate plotted in the middle panels. The nutrition rates are identical in all subfigures since we do not change nutrition rate. Panels \textbf{(A)} and \textbf{(B)} represent a virtual patient’s response to LQG-suggested and protocol-suggested IV insulin rates. Panel \textbf{(C)} shows a possible GM result if we had a BG model that could represent the real BG dynamics of ICU patients perfectly.}
	\label{simulation}
\end{subfigure}

We included nutrition change times as intervention times for the LQG controller in addition to the intervention times of a protocol. We see the positive effect of this strategy by comparing the BG values between 15-20 hours. Since the protocol does not account for nutrition amount changes, it continually suggested higher IV insulin rates while the LQG controller suggested lower rates due to reduced nutrition amounts over that time interval. This difference shows the importance of accounting for nutritional changes in developing BG control strategies.

Figure \ref{lqg_msgm_sim} shows simulated BG values by the MSG model in response to LQG-suggested IV insulin rate. This figure demonstrates potential GM performance given a perfect representation of BG dynamics with a physiology-based mechanistic model. Our target value for the LQG controller is the upper limit of the respective protocol's target range (180 mg/dl in Figure \ref{simulation}) to avoid hypoglycemia. The BG values simulated with the MSG model hover around that target value, with occasional lower values due to nutrition rate changes shifting the system equilibrium. However, the LQG controller quickly identifies the IV insulin rate to shift the equilibrium to the target value of 180 mg/dl. Therefore, one way to develop an effective model-based BG controller is to use a model that accurately represents BG dynamics. However, such models do not exist and might not be useful in control frameworks because of their complexity and real-world data limitations. Hence, we aim to develop controllers using simpler models accounting for the system's dynamic behavior and real-world data limitations through novel techniques (such as introducing resolvable stochasticity into the model) rather than increasing model complexity.